\documentclass[12pt]{article}
\usepackage{amsmath,amsfonts, amssymb,verbatim}
\newtheorem{pkt}{}[section]
\newcommand{\bpk}{\begin{pkt}\rm }
\newcommand{\epk}{\end{pkt}}

\newcommand{\M}{\mathbb{M}}

\newcommand{\R}{{\mathbb R}}

\newcommand{\N}{{\mathbb N}}
\newcommand{\C}{{\mathbb C}}

\newcommand{\A}{\mathbb{A}}
\newcommand{\B}{\mathcal{B}}
\newcommand{\be}{\begin{equation}}
\newcommand{\ee}{\end{equation}}

\newcommand{\HH}{\mathcal{H}}

\newcommand{\ra}{\rangle}
\newcommand{\la}{\langle}

\begin{document}
\title{Axioms of Quantum Mechanics  in light of Continuous Model Theory }
\author{B.Zilber}
\maketitle
%\tableofcontents

\begin{abstract}
We revisit Dirac's  axiomatization of quantum mechanics and show that it can essentially be   reformulated in a more familiar  logical setting - continuous model theory. Our aim is twofold: (i) to present the Dirac--von Neumann formalism in a genuinely axiomatic manner suitable for logicians, and (ii) to exhibit a structural analogy between Hilbert spaces and Tarski's cylindric algebras, which were introduced in the program of algebraisation of first-order logic. 

Recall that the cylindric algebra $\mathfrak{C}(\A)$ of a first order structure $\A$ allows to recover $\A$ up to elementary equivalence.

For a general continuous structure $\M$, we introduce an analogue $\B(\M)$ of the cylindric algebra  of a first--order structure. Under natural tameness assumptions, $\B(\M)$ takes the form of a (rigged) Hilbert space with operators, and $\M$ can be recovered from $\B(\M)$.

\end{abstract}

\section{Introduction}
\bpk 

The axiomatic formulation of quantum mechanics was introduced by Paul Dirac in 1930 \cite{Dirac} through a description of a Hilbert space, and later developed with greater mathematical rigor in a monograph of 1932 by John von Neumann. Since 1930, Dirac went through several rewritings and new editions to refine his calculus to a level he considered satisfactory. In the 1950s the theory settled with the notion of {\em rigged Hilbert space}, or Gel'fand triple. 
 Modern books present Dirac's axioms in a succinct form, often omitting much of the technical detail.

In section \ref{s2}, for the benefit of the reader, we take a glance at   the  axioms of quantum mechanics  for physicists as presented in  \cite{ModernP}. 
Readers with a background in logic  will notice that what physicists call ``axioms'' is very far from what is a conventional set of axioms  in a formal language, even in its early form such as Hilbert's axiomatisation of geometry \cite{HilbertGrundlagen}. 

In section \ref{s4} we argue  that the language that Dirac informally used is that of {\em continuous logic}. We then go further and explain that  Dirac's axiomatisation has chosen the formalism termed   {\em algebraic logic} as exemplified e.g. by  A.Tarski's cylindric algebras - the tower of Lindenbaum algebras $\mathcal{F}_n(\A)$ of a first order structure $\A$ expanded by quantifiers,  \cite{Monk}. The cylindric algebra of $\A$ characterises $\A$ uniquely up to elementary equivalence.
%In fact, Hilbert spaces can be seen as a continuous model theory version of cylindric algebras. 

The main contribution of the current paper is the introduction and an initial study
 of an analogue $\B(\M)$ of a cylindric algebra  
for a structure $\M$ of continuous logic, thereby showing that the algebraisation of continuous logic naturally leads to a Dirac - von Neumann style formalism. 

In general,  $\B(\M)$ is a tower of Banach spaces $B_n(\M)$ with linear operators -  an analogue of the tower of Lindenbaum algebras. 
Theorem \ref{MainT} proves that under stronger assumptions on $\M$, $\B(\M)$ acquires the structure of  a tower of (rigged)  pre-Hilbert spaces expanded by linear operators. Crucially, in this case $\M$ can be recovered canonically from $\B(\M).$

Completing the spaces in $\B(\M)$ with respect to their inner product norms yields a corresponding tower of Hilbert spaces associated with $\M.$
Hilbert spaces are, of course, a central object of study in continuous model theory. Yet it is worth noting that, for example, if 
$\M$
is a continuous structure on a compact universe (corresponding to the ``particle in a box'' model in quantum mechanics), then 
$\M$
is the unique model of its theory -  a case of absolute categoricity. In contrast, the theory of the associated Hilbert space expanded by non‑commuting operators, is known to be highly complicated.
\epk
\bpk Continuous logic and continuous model theory were introduced in the monograph \cite{CK} in the 1960s and have since  been developed and further generalised for various applications. For readers without a background in logic  the article of E.Hushovski \cite{HrItamar} outlines  a philosophy behind the mathematical formalism. 

The link between physics formalism and continuous logic was proposed and initially explored by the present author in \cite{ZCL}. Paper \cite{DvN} is devoted to a rigorous interpretation of Dirac - von Neumann axioms as a system of  axioms of continuous logic and the study of its models. 
\epk

\section{Dirac's calculus and axiomatisation of\\  quantum mechanics}\label{s2}

In order to make a link to the definitions and statements in the rest of the paper, 
 we reproduce a slightly edited version of axioms as presented for physicists in \cite{ModernP}, 6.3.

\bpk
{\bf Axiom 1}. The “state” of a quantum system is described by a non-zero vector $|\psi\ra$ belonging to
a complex Hilbert space $\HH.$ This state is usually called ``ket $\psi$''.
A complex Hilbert space $\HH$ is a vector space, which can be finite dimensional
or infinite dimensional, equipped with the complex scalar product (also called inner
product) $\la \psi |\psi' \ra$ between any pair of states $|\psi\ra,$ $|\psi'\ra$ in $\HH.$ The norm, or modulus, of
a generic vector $|\psi\ra\in \HH$ is defined as
$$||\psi||=|\la\psi  |\psi\ra|$$
and usually $|\psi\ra$ is normalized to one, i.e.$||\psi||=1.$ The symbol $\la \psi|$ which appears in
the definition of the norm is called ``bra $\psi$'' and it can be interpreted as the function $$\la \psi|: \HH\to \C.$$
 For any $|\psi'\ra\in \HH$  this function gives a complex number $\la \psi|\psi'\ra$ obtained as
scalar product of $|\psi\ra$ and $|\psi'\ra$. In a complex Hilbert space $\HH$ there exists a set of basis
vectors $|\phi_\alpha\ra$  which are orthonormal, 
i.e. $\la \phi_\alpha|\phi_\beta\ra =\delta(\alpha-\beta),$ and such that \be \label{psi}|\psi\ra=\sum_\alpha c_\alpha |\phi_\alpha\ra\ee
for any $|\psi\ra$, where the coefficients $c_\alpha$ belong to $\C.$

\medskip

{\bf Axiom 2}. Any observable (measurable quantity) of a quantum system is described
by a self-adjoint linear operator $F: \HH\to \HH$ acting on the Hilbert space of state vectors.

For any classical observable $F$ it exists a corresponding quantum observable $F$.

\medskip

{\bf Axioms 3}. The possible measurable values of an observable $F$ are its eigenvalues $f,$
such that
$$ F|f\ra=f |f\ra$$
with $|f\ra$ the corresponding eigenstate.
The observable $|f\ra$ admits the spectral resolution
\be \label{f}F=\sum_f f |f\ra \la f|\ee 
where  $\{ |f\ra\}$ is the set of orthonormal eigenstates of $F$, and the mathematical object $\la f|$, called ``bra of $f$'', is a linear map that maps into the complex number. This
also satisfy
 the identity
$$\sum_f |f\ra \la f|=\mathrm{I}.$$

\medskip

{\bf Axiom 4}. The probability $P$ of finding the state $|\psi\ra$ in the state $|f\ra$ (both of norm 1) is given by
$$P = |\la f|\psi\ra|^2$$
 This probability $P$ is also the probability of measuring the value $f$ of
the observable $F$ when the system is in the quantum state  $|\psi\ra$.

\medskip

{\bf Axiom 5}. The time evolution of states and observables of a quantum system with
Hamiltonian $\mathrm{H}$ is determined by the unitary operator
$$K^t:=\exp(-i\mathrm{H}t/\hbar)$$   
such that $|\psi(t)\ra  = K^t|\psi\ra$ 
is the time-evolved state $|\psi\ra.$ 
\epk
\bpk  Now we make several comments on the axioms.

The term ``Hilbert space'' here should actually be read  as {\em the rigged Hilbert space} (see \cite{rigged}), equivalently, a {\em Gel'fand triple}.  It differs from the standard definition by accommodating along with a Hilbert space $\HH$ both a subspace $\Phi$ of {\em test-functions}  and the dual space $\Phi^*$ with 
$$\Phi\subset \HH\subset \Phi^*.$$

%The position states $|x\ra$ form the orthonormal basis of the space  $\Phi^*$ of bra-vectors.

The summation formulas like (\ref{psi}) and (\ref{f}) are presented in a form of an integral if the family $|\psi_\alpha\ra$  is continuous but seems natural in the summation form when  $\alpha$ runs in the discrete spectrum of an operator.

\epk
\bpk\label{R3.4} {\bf Remark.} 
Rigged Hilbert spaces provide a powerful mathematical framework to extend quantum mechanics, allowing distributions and generalized eigenfunctions to be rigorously handled. However, it is still contested whether the condition of completeness of $\HH$ is necessary. It  is almost generally accepted that not every element corresponds to a physically realisable state -- some are purely mathematical artifacts, see e.g. \cite{Nonphysical}. 

In the more general context of quantum field theories Wightman axioms explicitly postulate that physically meaningful part of the rigged Hilbert space $\HH$ is a dense subset $\mathcal{D}\subset \HH.$ 
\epk

\section{Algebraisation of Logic}\label{s4}

We recall the algebraisation of first--order logic via cylindric algebras,  see \cite{Monk}.

\bpk {\bf Cylindric algebras}

Let $\A$ be a first--order structure in a language $\mathcal{L}$. For each finite tuple of distinct variables $\bar{x}=(x_{i_1},\dots,x_{i_n})$, let $\mathcal{F}_{\bar{x}}$ be the Lindenbaum--Tarski algebra of $\mathcal{L}$--formulas in variables $\bar{x}$ modulo equivalence in $\A$.

There are natural embeddings $\mathcal{F}_{\bar{x}}\hookrightarrow\mathcal{F}_{\bar{x}'}$ whenever $\bar{x}\subseteq\bar{x}'$. Define $$ \mathcal{F}:=\bigcup_{\bar{x}} \mathcal{F}_{\bar{x}}. $$

\paragraph{\textbf{Cylindrification.}} For each variable $x_k$, define $$ \exists x_k:\mathcal{F}_{\bar{x}{x_k}}\to\mathcal{F}_{\bar{x}} $$ by $$ \exists x_k[\varphi(\bar{x},x_k)] := \text{equivalence class of }\exists x_k,\varphi(\bar{x},x_k). $$

The \emph{cylindric algebra} of $\A$ is $$\mathfrak{C}(\A):=(\mathcal{F},{\exists x_k}_{k\in\mathbb{N}}). $$ 
\epk 
\bpk {\bf Representation theorem} \label{TMH}
 [Henkin--Monk--Tarski]
 {\em 
  Let $\A$ and $\A'$ be structures in the same first--order language. Then $$ \A\equiv \A' \quad\Longleftrightarrow\quad \mathfrak{C}(\A)\cong\mathfrak{C}(\A'),$$ where the isomorphism identifies  sets definable by the same formulas.}
\epk

\section{Algebraisation of Continuous Logic}

\bpk{\bf Continuous language and structures}

The definitions are basically the same as in \cite{BBHU} except for the fact that the range of predicates are bounded domains in $\C:$ 
 
{\bf A continuous language $\mathcal{L}$} consists of: \begin{itemize} \item a distinguished binary predicate $\mathrm{dist}(x,y)$ with values in a closed interval of $\R,$ \item predicate symbols $P$ interpreted in a model $\M$ with universe $M,$ a metric space, as uniformly continuous maps $P^M:M^n\to D_P$, a compact subset of $\C$, \item function symbols $f$ interpreted as uniformly continuous maps $f^M:M^m\to X$, \item continuous connectives $\alpha:D_1\times\cdots\times D_k\to D$, \item quantifiers such as $\sup$ and $\inf$ but can be others which are expressible in terms of the basic ones. \end{itemize}

{\bf A continuous $\mathcal{L}$--structure} $\M$ consists of: \begin{itemize} \item a finite diameter complete metric space $(M,\mathrm{dist})$, \item interpretations of predicate and function symbols as above.
\item we also consider a more general setting with countably many universes/ metric spaces $$M_1\subseteq M_2\subseteq \ldots\subseteq M_n\subseteq \ldots$$ 
 of finite diameters, such that $\mathrm{dist}$ on $M_{i}$ extending    $\mathrm{dist}$ on $M_{i+1}$.

 \end{itemize}

{\bf Formulas, connectives, quantifiers}

Formulas are built from atomic predicates  using substitutions of terms, connectives and quantifiers.
% Quantifiers are uniformly continuous operators $Q:\mathcal{F}_{n+1}\to\mathcal{F}_n$.
\epk

\bpk {\bf Definable predicates, maps, sets and imaginaries}\label{4.2}

A {\bf definable predicate} on $M$ is a uniformly continuous map $\varphi:M^n\to D_\varphi$ obtained as a uniform limit of interpretations of formulas.

A set $A\subseteq M^n$ is definable if $A=\{ a:\varphi(a)=0\}$ for a definable predicate $\varphi$ satisfying the almost--near condition. 

A map $f:M^m\to M^n$ is definable if its graph is definable. 

\medspace

Recall  a continuous structure $\M$ with universe $M$ can be extended to the continuous multisorted structure $\M^\mathrm{eq}$ with (imaginary) sorts of the form $A/E,$ $A\subset M^n$ a definable set and $E(x,y)$ a binary relation on $A$ which is a pseudo-metric. In particular, given an $m+n$-ary predicate $\psi(x,y)$ one defines 
$$E_\psi(y_1,y_2)=\sup_x |\psi(x,y_1)-\psi(x,y_2)|$$
Then $M^n/E_\psi$ is the set of imaginary elements called canonical parameters for $\psi(x,y)$ and $$\{ \psi(x, c): \ c\in M^n/E_\psi \} $$ is a {\bf definable family of definable sets} which is in bijective correspondence with the set of canonical parameters $M^n/E_\psi$ and each $c$ a point in $\M^\mathrm{eq}.$  

In particular, if $n=0$ and $\psi(x)$ is a definable $m$-ary predicate, there is a unique point $c_\psi\in \M^\mathrm{eq}$, a canonical parameter for $\psi.$  
%The set $  B_n(\M)$ of canonical parameters   
\epk 
     
\bpk {\bf Spaces of definable predicates}

Let $  B_n(\M),$ $n\ge 1,$ be the set of definable predicates $\varphi:M^n\to D_\varphi.$ Set $B_0(\M)=\C.$ 
By the last remark in \ref{4.2} we may treat $  B_n(\M)$ as a subset of $\M^\mathrm{eq}$ with each point marked by a constant symbol $c_\varphi$ and  with the metric of the normed space,  the norm defined by $$|\varphi|:=\sup_{x\in M^n}|\varphi(x)|.$$

 Given a definable map $$f: M^n\to M^m, \  (m,n\ge 0)$$ define the map between the spaces of predicates
 $$L_f:   B_m(\M)\to   B_n(\M); \ \   \psi \mapsto \psi\circ f$$%  \ (\psi\circ f)(\bar{x})=\psi( f(\bar{x}))$$
 
 Given a predicate $\phi\in   B_n(\M)$  define the map between the spaces of predicates
 $$L_\phi:   B_n(\M)\to   B_n(\M); \ \   \psi \mapsto \phi\cdot\psi$$%
\epk
\bpk 
{\bf Proposition.} {\em  For each $n,$ $  B_n(\M)$ is a complex Banach space and $L_\phi,$  
	$L_f$ are continuous linear operators between the spaces. 

Suppose $\M'$ is a structure in the same language and $\M'\equiv \M$ as continuos structures. Then there is a unique isomorphism
 $$  B_n(\M)\cong   B_n(\M')$$ of Banach spaces with operators.  }

\medspace

{\bf Proof.} By the general continuous model theory definable predicates are uniformly continuous and have bounded codomain. Hence
finite linear combinations of definable predicates are definable and have a finite norm and a Cauchy sequence of definable predicates has a limit which  is uniformly continuous and definable. 

The statement on $L_\phi$ and $L_f$ is a direct consequence of assumptions.

Finally, note that the map $$c_\phi^{\M} \mapsto c_\phi^{\M'}$$ is an   isomorphism of Banach spaces which also preserves the linear operators. 
$\Box$

\medspace 

{\bf Remark 1.}  The structure on $  B_n(\M)$ also encodes the multiplicative structure. For $\varphi,\psi\in   B_n(\M):$ $$ \varphi\cdot \psi= L_\varphi(\psi).$$   
%This distance coincides with the distance on the respective set of points in the continuous structure $\M^\mathrm{eq},$ %namely, for $\psi, \varphi\in   B_n(\M)$
%$$\mathrm{dist}(c_\psi, c_\varphi)=\sup_x |\psi(x) - \varphi(x)|.$$

\medspace

{\bf Remark 2.}  Consider the tower of Banach spaces with operators
\be \label{Btower} \mathcal{B}(\M):= \{ B_1(\M)\hookrightarrow  B_2(\M)  \hookrightarrow \ldots \hookrightarrow   B_n(\M) \ldots   \}\ee 
where $  B_n(\M)\hookrightarrow   B_{n+1}(\M)$ are all possible  embeddings identifying $\psi \in   B_n(\M)$ as an $n+1$-variable predicate in $   B_{n+1}(\M)$ with an extra fictitious variable.

 The proposition above shows that $\mathcal{B}(\M)$
%of the form $\psi(\bar{x})\mapsto \psi'(\bar{x'})  
is a good candidate for a ``cylindric algebra of a CL-structure $\M$''.   Banach spaces are well-studied objects of continuous model theory and we will treat $\mathcal{B}(\M)$ as a multisorted continuous  structure with operators -  definable maps. 

However, this construction is not sufficient to provide a good analogue of the Henkin - Monk - Tarski Theorem \ref{TMH}. One needs a  duality theory which is available in the setting of Hilbert spaces.  

\epk 

\bpk\label{3.7} {\bf Continuous linear functionals on $\  B_n(\M)$. }

A continuous linear map $$I:   B_n(\M) \to \C$$ is called a linear functional on $  B_n(\M).$

\medspace

{\bf An evaluation functional} on $  B_n(\M)$ is a map $  B_n(\M)\to \C$ of the form
$$\lambda_x: \psi \to \psi(x),\mbox{  for } x \in \M^n.$$

It is immediate from the definitions that an evaluation functional is a continuous linear functional. 

The linear space of continuous linear functional on  $  B_n(\M)$ will be called  {\bf the dual space} $ B^*_n(\M).$
%Let $\Phi\subset B(X)$ be a definable family of predicates. %
%For each $x\in X$, define $\lambda_{x,\Phi}:\Phi\to\mathbb{C}$ by $\lambda_{x,\Phi}(\varphi)=\varphi(x)$.

\medspace

Recall the following classical facts:

\medspace

{\bf Fact 1}. (Riesz representation theorem.) {\em Every continuous linear functional $I$ on a space of continuous functions comes from integration against a measure:  
	\be \label{L}I: \psi \mapsto \int \psi d\mu_I,\ \ \mu_I\in \mathcal{M}(M^n) \ee
	where  $\mathcal{M}(M^n)$ is the space of complex (or finite signed) regular Borel measures on $M^n.$ 
}

\medspace

{\bf Fact 2}. %Regular Borel measure $\mu$ on a compact metric space $\Omega$ is weakly approximated by atomic 
Finite atomic measures are weak-$^*$
dense in $\mathcal{M}(M)$ for compact $M.$

In particular, given $\mu,$ for any system of $\frac{1}{N}$-dense lattices $M_N\subset M,$ $N\in \N,$ there is a system $\mu_N$ of  atomic measures based on points of $M_N$ such that 
for any continuous $\psi$ on $M$ 	\be \label{sunmu}\int \psi d\mu= \lim_{N\to \infty} \sum_{x\in M_N} \psi(x) \mu_N(x)\ee 

\medspace

{\bf Remark.} In quantum mechanics evaluation functionals are presented as the Dirac delta-functions $\delta_x.$

\epk
\bpk Let $M_N$ be an $\frac{1}{N}$-dense lattice on $M$ and let 
the atomic measure $\nu_N$ given by
 $$\nu_N(x)=\left\lbrace\begin{array}{ll}  \frac{1}{N}, \mbox{ if } x\in M_N\\ 
0, \mbox{ otherwise }\end{array}\right.  
$$ 
%will be called canonical for the lattice $M_N.$ 

The measure $\nu$ defined as the limit of $\nu_N$ by formula (\ref{sunmu}) on a compact $M$ will be called {\bf  the canonical measure}  on  $M.$ 

\medspace

{\bf Remark.} Assume that for each number $N\in \N$ there is a set $M_N\subset M$ of definable points  which is $\frac{1}{N}$-dense in $M.$ Then every point of the universe $M$  of the continuous structure $\M$ is definable. Indeed, definable points of $\M$ are dense in $M$ and a limit of a sequence of definable points is definable.

%Under this assumption integration over the canonical measure is definable. 

\epk
\bpk The next theorem is formulated under assumptions on $\M$ that might look too restrictive for continuous model theory. However, in quantum mechanics the standard assumption is that the Hilbert space in question is the space of wave functions on a real manifold $M$, even more specifically $M=\R^3,$ and in the so called 1-dimensional quantum mechanics, $M=\R.$  This is the setting of unbounded continuous model theory which is readily reducible to the case of bounded many sorted continuous structures by choosing the sorts $$[-1,1]\subset [-2,2]\subset\ldots [-n,n]\subset \ldots $$

Note that in such a structure, assuming the end-points of $M=[-1,1]$ are definable, every point $x$ is definable by formulas 
 $\mathrm{dist}(-1,x)=d_x^-$ and $\mathrm{dist}(1,x)=d_x^+,$ for respective $d_x^-,d_x^+\in \R_{\ge 0}.$

\medspace

We say that a continuous structure $\M$ is {\bf tame} if the universe $M$ is compact and all points are definable.
\epk

\bpk \label{MainT}
{\bf Theorem}. {\em Let $\M$ be a tame continuous structure.

		(A).   	Any continuous linear functional   on $  B_n(\M)$  is defined 
		over $\M$ by formula (\ref{sunmu}) for an appropriate measure $\mu.$ 
		
	 The 	integration operator $$\psi \mapsto \int \psi d \nu; \   B_{n+1}(\M)\to   B_n(\M)$$ over the canonical measure on $M$ is defined over $\M.$ 
	 
	\medspace
	
	(B). 
	The  set of evaluation  functionals on $  B_n(\M)$
	$$(M^n)^*:=\{ \lambda_{\bar{a}}: \ \bar{a}\in M^n\}$$ is in bijective correspondence with
	$M^n$
	%There is a bijection 
	$$\pi_n: (M^n)^*\to M^n; \ \ \lambda_{\bar{a}} \mapsto  \bar{a}.$$
	 
	 The continuous structure $\M$ is determined uniquely, up to isomorphism, by  $ \B(\M)$ together with evaluation functionals on $ B_n(\M),$ $n\in \N.$

	\medspace

	(C). 
	There is a canonical embedding, for each $n\in \N,$ 
	$$i_n:  B_n(\M) \hookrightarrow  B^*_n(\M).$$ 
	%via $$\alpha\mapsto \phi_\alpha,\ \ \alpha \in \B(\Omega^n), \ \phi_\alpha\in \B(\Omega^n)^*,
	%\ \phi_\alpha: \psi\mapsto \int \psi\cdot \bar{\alpha}\, d\nu.$$ 
	
	 $  B_n(\M)$ equipped  with the Hermitian inner product $$\la \psi|\alpha\ra:= \int \psi\cdot \bar{\alpha}\, d\nu^n$$ 
	 is a pre-Hilbert space 
	and its completion is a  Hilbert space  $L^2(M,\nu) $ embedded naturally into    $ B^*_n(\M):$ 
	$$B_n(\M) \subset L^2(M^n,\nu^n) \subset B^*_n(\M).$$ 
	
	The bijection $\pi=\{ \pi_n: n\in \N\}$ of (B)  induces an isomorphism from the structure $\M$ to a structure $\M^*$ based  on universe $M^*$ and predicates defined in terms of $ B_n^*(\M).$
	}
\medspace

{\bf Proof.} (A). This is just a summary of Facts 1 and 2 of \ref{3.7}.
%The first claim is by Fact 2, formula (\ref{sunmu}) using the assumption that $M$ is compact and every point is definable.

%The second claim is a corollary of the first.

\medspace

(B). First note that %a point $\psi\in  B_n(\M)$ induces a  predicate on $(M^n)^*:$ $$\psi: \lambda_{\bar{a}} \mapsto \lambda_{\bar{a}}(\psi); \ \  (M^n)^*\to \C.$$
%The 
predicates from $ B_1(\M)$ separate points of $M,$ that is, for each distinct $a,b\in M$ there is $\phi\in  B_1(\M)$ such that $\phi(a)\neq \phi(b).$ It is enough to take $\phi(x):= \mathrm{dist}(a,x).$ 

It follows that  $n$-predicates of the form   $\phi_1(x_1)\cdot \ldots\cdot\phi_n(x_n),$
$\phi_1,\ldots,\phi_n\in  B_1(\M)$, separate points of $M^n.$ 
 
Thus, by  the Stone - Weierstrass theorem, the linear subspace of $ B_n(\M)$ spanned by $n$-predicates of the form   $\phi_1(x_1)\cdot \ldots\cdot\phi_n(x_n)$ is dense in $ B_n(\M).$ 

Since for any $\phi_1,\ldots,\phi_n\in  B_1(\M),$ for each $\la a_1,\ldots,a_n\ra\in M^n,$
$$ \lambda_{a_1\ldots a_n} [  \phi_1(x_1)\cdot \ldots\cdot\phi_n(x_n)] = \lambda_{a_1}[\phi_1(x)]\cdot\ldots\cdot  \lambda_{a_n}[\phi_n(x)]$$
the linear functional  $\lambda_{a_1\ldots a_n}$ on $ B_n(\M)$ is determined by the $n$-tuple of linear functionals $\lambda_{a_1},\ldots, \lambda_{a_n}$ on $ B_1(\M).$

Thus, we proved that the correspondence $\lambda_{\bar{a}}\mapsto \bar{a}$ induces bijections 
$$(M^n)^* \to M^n, \ (M^n)^*\to (M^*)^n.$$

Finally, we define: 

the distance predicate on $M^*$ by
$$\mathrm{dist}(\lambda_{a_1}, \lambda_{a_2}):= \lambda_{a_1,a_2}[\mathrm{dist}(x_1,x_2]).$$

the $n$-ary predicate $\varphi$ on $M^*$ by
$$\varphi(\lambda_{a_1},\ldots, \lambda_{a_n}):= \lambda_{a_1,\ldots,a_n}[\varphi(x_1,\ldots,x_n]).$$

the definable map on $(M^*)^m\to (M^*)^n$ by
$$f(\lambda_{a_1},\ldots, \lambda_{a_m}):= \lambda_{a_1,\ldots,a_m}\circ L_f$$
for all $a_1,\ldots,a_n,\ldots \in M.$

By definition we then get
$$\mathrm{dist}(\lambda_{a_1}, \lambda_{a_2})=\mathrm{dist}(a_1,a_2]).$$
$$\varphi(\lambda_{a_1},\ldots, \lambda_{a_n})= \varphi(a_1,\ldots,a_n])$$	
and 
$$f(\lambda_{a_1},\ldots, \lambda_{a_m})= f{a_1,\ldots,a_m}$$
which proves that $$\pi:  M^*\to M$$ is an isomorphism of the continuous structures.

\medspace

(C).  
 The formula for $\la \psi|\alpha\ra$ satisfies the assumptions of Hermitian inner product on $\B_n(\M)$ by definitions, thus $B_n(\M)$ is a pre-Hilbert space.
 
   By the standard argument, for each $\alpha\in B_n(\M)$ the map
   $$I_\alpha: \psi \mapsto \la \psi|\alpha\ra; \ B_n(\M)\to \C$$
   is a continuous linear functional, $I_\alpha\in B_n^*(\M).$ Thus $\alpha \mapsto I_\alpha$ is a map $i_n: B_n(\M)\to B_n^*(\M).$ This map is an embedding since $I_{\alpha}=0$ if and only if $\alpha=0.$ 

Now note that the completion of $B_n(\M)$ in the inner product topology is the Hilbert space of square integrable functions over $M^n,$ that is the space $L^2(M^n,\nu^n).$ It is a well-known fact that the Dirac delta functions $delta_a$ (that is $\lambda_a$) are in $B_n^*(\M)$ but not in $L^2(M^n,\nu^n).$

$\Box$

\epk

\bpk {\bf Remarks.}  For a tame $\M$  we get $n$-Gel'fand triples
$$B_n(\M)\subset L^2(M^n, \nu^n) \subset B_n^*(\M)$$
 
This system can be taken to be the cylindric algebra  of the continuous structure $\M.$

\epk

\end{document}